\documentclass[a4paper,11pt]{article}
\usepackage{graphics}
\usepackage{fullpage}
\usepackage{amsmath}
\usepackage{amssymb}
\usepackage{theorem}

\DeclareMathOperator{\Rot}{Rot}

\DeclareMathOperator{\InfX}{\mbox{\upshape\sffamily\bfseries T}}
\newcommand{\InfG}{\InfX^{\circ}}

\newcommand{\sigmap}[1][T]{\ensuremath{\mathcal{C}_1^{\sigma}(#1)}}
\newcommand{\TR}{T_{\IR}}
\newcommand{\RotR}{\Rot_{\IR}}

\newcommand{\rhos}[1][F]{\rho_{_{#1}}}
\newcommand{\orhos}[1][F]{\overline{\rho}_{_{#1}}}
\newcommand{\urhos}[1][F]{\underline{\rho}_{_{#1}}}

\newcommand{\pluscover}[1]{\displaystyle \mathop{\longrightarrow}_{#1}^{+}}

\newcommand{\modi}{({\rm mod\ 1})}
\newcommand{\eps}{\varepsilon}
\newcommand{\Int}[1]{{\rm Int}\left(#1\right)}
\newcommand{\diam}[1]{\mbox{\rm diam}\left(#1\right)}

\newcommand{\IN}{\ensuremath{\mathbb{N}}}
\newcommand{\IZ}{\ensuremath{\mathbb{Z}}}
\newcommand{\IQ}{\ensuremath{\mathbb{Q}}}
\newcommand{\IR}{\ensuremath{\mathbb{R}}}

\newcommand{\CC}{{\cal C}}
\newcommand{\CD}{{\cal D}}

\newtheorem{theo}{Theorem}[section]
\newtheorem{prop}[theo]{Proposition}
\newtheorem{lem}[theo]{Lemma}

\theorembodyfont{\rm} 
\newtheorem{defi}[theo]{Definition}
\newtheorem{ex}[theo]{Example}
\newtheorem{rem}[theo]{Remark}

\newenvironment{demo}{\medskip\noindent {\it Proof}.}{\hfill$\Box$\par\medskip}

\begin{document}

\title{Rotation set for maps of degree 1 on the graph sigma}
\author{Sylvie Ruette}
\date{December 7, 2007}


\vspace*{1mm}

\begin{center}
{\LARGE Rotation set for maps of degree 1 on the graph sigma}

\bigskip
{\large Sylvie Ruette}

\medskip
December 7, 2007
\end{center}

\begin{abstract}
For a continuous map on a topological graph containing a unique loop
$S$ it is possible to define the degree and,  for a map of degree $1$,
rotation numbers.  It is known that  the set of rotation numbers of
points in $S$ is a compact interval  and for every rational $r$
in this interval there exists a periodic point of rotation number
$r$.  The whole rotation set (i.e. the set of all rotation
numbers) may not be connected and it is not known in general whether
it is closed.

The graph sigma is the space consisting in  an interval attached by
one of its endpoints to a circle.  We show that, for a map of degree 1
on the graph sigma,  the rotation set is closed and has finitely many
connected components. Moreover, for all  rational numbers $r$ in the
rotation set, there exists a  periodic point of rotation number $r$.
\end{abstract}

\section{Introduction}

\begin{figure}[htb]
\centerline{\includegraphics{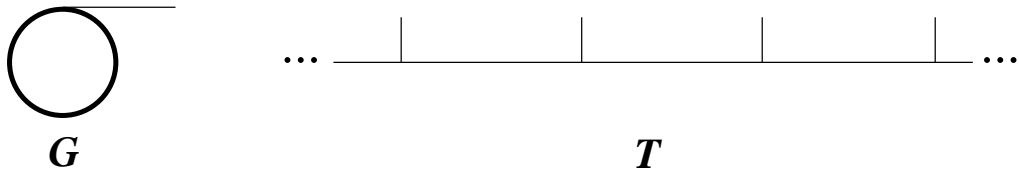}}
\caption{$G$ is the graph $\sigma$, its universal covering is $T$.
\label{fig:sigma}}
\end{figure}

\begin{figure}[hbt]
\centerline{\includegraphics{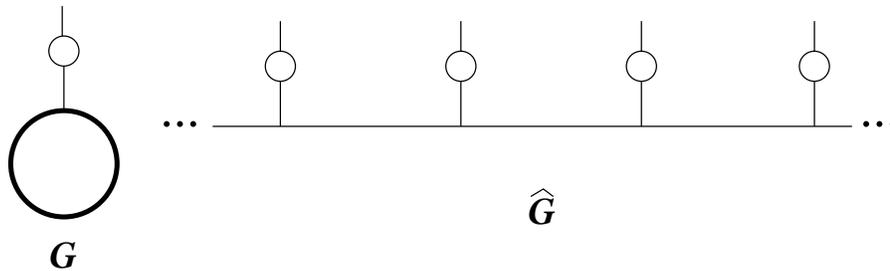}}
\caption{The graph $G$ is unwound with respect to the bold loop
to obtain $\widehat{G}$, which is a lifted graph. \label{fig:hatG}}
\end{figure}

In \cite{AR} a rotation theory is developed for continuous self maps
of degree 1 of topological graphs having a unique loop, using the ideas
and techniques of \cite{Mis, BGMY}.
A rotation theory is usually developed in the universal covering space by using
the liftings of the maps under consideration. The universal covering of a graph
containing a unique loop is an ``infinite tree modulo 1'' (see
Figure~\ref{fig:sigma}). It turns out that the rotation
theory on the universal covering of a graph with a unique loop can be easily
extended to the setting of infinite graphs 
that look like the space $\widehat{G}$
from Figure~\ref{fig:hatG}. These spaces are defined in detail in
Section~\ref{subsec:liftedgraphs} 
and called \emph{lifted graphs}. Each lifted graph
$T$ has a subset $\widehat{T}$ homeomorphic to the real line $\IR$ that
corresponds to an ``unwinding'' of a distinguished loop of the original
graph. In the sequel, we identify $\widehat{T}$ with $\IR$.

Given a lifted graph $T$ and a map $F$ from $T$ to itself of degree one, there
is no difficulty to extend the definition of rotation number to this setting in
such a way that every periodic point still has a rational rotation number as in
the circle case. However, the obtained rotation set $\Rot(F)$ may not be
connected. Despite of this fact, it is proven in \cite{AR} that the
set $\Rot_{\IR}(F)$ corresponding to the rotation numbers of all points
belonging to $\IR$, has properties which are similar to (although weaker than)
those of the rotation interval for a circle map of degree one. Indeed, this set
is a compact non empty interval,
if $p/q\in \Rot_{\IR}(F)$ then there exists a periodic point of rotation number
$p/q$, and if $p/q\in \Int{\Rot_{\IR}(F)}$ then for all large enough positive
integers $n$ there exists a periodic point of period $nq$ of rotation number
$p/q$.  

We conjecture that the whole rotation set $\Rot(F)$ is closed.
In this paper, we prove that, when the
space $T$ is the universal covering of the graph $\sigma$ consisting in 
an interval attached by one of its endpoint to a circle
(see Figure~\ref{fig:sigma}), then the rotation set is the union of 
finitely many compact intervals. Moreover, all 
rational points $r$ in $\Rot(F)$ are rotation numbers of periodic
{\modi} points. It turns out
that the proofs extend to a class of maps on graphs that we call
$\sigma$-like maps, which are defined in Section~\ref{subsec:sigmalikemaps}.

\section{Definitions and elementary properties}
\subsection{Lifted graphs}\label{subsec:liftedgraphs}
A \emph{(topological) finite graph} is a compact connected set $G$ containing a
finite subset $V$ such that each connected component of $G\setminus V$
is homeomorphic to an open interval. 

The aim of this section is to define in detail the class of \emph{lifted 
graphs} where we develop the rotation theory. They are obtained from 
a topological graph by unwinding one of its loops. This gives a new 
space that contains a subset homeomorphic to the real line and that is 
``invariant by a translation'' (see Figures~\ref{fig:sigma} and 
\ref{fig:hatG}). In \cite{AR}, a larger class of spaces called \emph{lifted
spaces} is defined.

\begin{defi}
Let $T$ be a connected closed topological space. We say that $T$ is a
\emph{lifted graph} if there 
exist a homeomorphism $h$ from $\IR$ into $T$,
and a homeomorphism $\tau\colon T\to T$ such that
\begin{enumerate}
\item $\tau(h(x)) = h(x+1)$ for all $x\in\IR$,
\item the closure of each connected component of  $T \setminus h(\IR)$ is a
finite graph that intersects $h(\IR)$ at a single point,
\item the number of connected components $C$ of $T \setminus h(\IR)$ such that
$\overline{C} \cap h([0,1])\neq \emptyset$ is finite.
\end{enumerate}
The class of all lifted graphs will be denoted by $\InfG$.
\end{defi}

To simplify the notation, in the rest of the paper we identify $h(\IR)$
with $\IR$ itself. In this setting, the map
$\tau$ can be interpreted as a translation by 1. So, for all $x \in T$ we
write $x+1$ to denote $\tau(x)$. Since $\tau$ is a homeomorphism, this notation
can be extended by denoting $\tau^m(x)$ by $x+m$ for all $m \in \IZ$.

Because of (ii), not all infinite graphs obtained by unwinding a finite graph 
with a distinguished loop are lifted graphs. The essential property of
this class is the existence of a natural retraction from $T$ to $\IR$.

\begin{defi}
Let $T\in \InfG$. The retraction $r\colon T\to\IR$ is the continuous map 
defined as follows. When $x\in\IR$, then $r(x) = x$. When
$x \notin \IR$, there exists a connected component $C$ of
$T \setminus \IR$ such that $x\in C$ and $\overline{C}$ intersects $\IR$ at a
single point $z$, and we let $r(x)=z$.
\end{defi}

\subsection{Maps of degree 1 and rotation numbers}

A standard approach to study the periodic points and orbits of a graph
map is to work at lifting level with the periodic {\modi} points.
The results on the lifted graph can obviously been pulled back to the original
graph (see \cite{AR}).
Moreover, the rotation numbers have a signification only for maps of degree 1,
as in the case of circle maps.
In this paper, we deal only with maps of degree 1 on lifted graphs.

\begin{defi}
Let $T \in \InfG$. A continuous map $F\colon T\to T$ is of degree $1$ if
$F(x+1)=F(x)+1$ for all $x\in T$.

A point $x\in T$ is called \emph{periodic {\modi}} for $F$ if there exists
a positive integer $n$ 
such that $F^n(x)\in x+\IZ$. The \emph{period} of $x$ is the
least integer $n$ satisfying this property, that is, $F^n(x)\in x+\IZ$
and $F^i(x)\not\in x+\IZ$ for all $1\leq i\leq n-1$.
\end{defi}

The next easy lemma summarises the basic properties of maps of degree $1$
(see for instance \cite[Section 3.1]{ALM}).

\begin{lem}\label{lem:degree1-Fn}\label{lem:compdeg1}
Let $T\in\InfG$ and $F\colon T\to T$ a continuous map of degree 1.
The following statements hold for $n\in\IN$, $k\in\IZ$
and $x\in T$:
\begin{enumerate}
\item $F^n(x+k)=F^n(x)+k$.
\item $(F+k)^n(x)=F^n(x)+kn$.
\item If $G\colon T\to T$ is another continuous map of degree 1, then $F\circ
G$ is a map of degree 1. In particular, $F^n$ is of degree 1 for all $n\geq 1$.
\end{enumerate}
\end{lem}

We define three types of rotation numbers.

\begin{defi}
Let $T\in\InfG$, $F\colon T\to T$ a continuous map of degree 1
and $x\in T$. We set
\[
\urhos(x) = \liminf_{n\to+\infty} \frac{r\circ F^n
(x)-r(x)}{n} \quad\mbox{and}\quad
\orhos(x) = \limsup_{n\to+\infty} \frac{r\circ F^n
(x)-r(x)}{n}.
\]
When $\urhos(x) = \orhos(x)$ then this number will be denoted by $\rhos{F}(x)$
and called the \emph{rotation number of $x$}. 
\end{defi}

We now give some elementary properties of rotation numbers
(see \cite[lemma 1.9]{AR}).

\begin{lem}\label{lem:first-properties}
Let $T\in\InfG$, $F\colon T\to T$ a continuous map of degree 1,
$x\in T$, $k\in\IZ$ and $n\in\IN$.
\begin{enumerate}
\item $\orhos(x+k)=\orhos(x)$.
\item $\orhos[(F+k)](x)=\orhos(x)+k$.
\item $\orhos[F^n](x)=n\orhos(x)$.
\end{enumerate}
The same statements hold with $\underline{\rho}$ instead of
$\overline{\rho}$.
\end{lem}

An important object that synthesises all the information about rotation numbers
is the \emph{rotation set} (i.e., the set of \emph{all} rotation numbers). Since
we have three types of rotation numbers, we have several kinds of rotation sets.

\begin{defi}
Let $T\in\InfG$ and $F\colon T\to T$ a continuous map of degree 1.
For  $S\subset T$ we define the following \emph{rotation sets}:
\begin{align*}
\Rot^+_S(F) &= \{\orhos(x)\mid x\in S\},\\
\Rot^-_S(F) &= \{\urhos(x)\mid x\in S\}, \text{ and}\\
\Rot_S(F)   &= \{\rhos(x)\mid x\in S\mbox{ and }\rhos(x)\mbox{ exists}\}.
\end{align*}
When $S=T$, we omit the subscript and we write $\Rot^+(F),\ \Rot^-(F)$ and
$\Rot(F)$ instead of
$\Rot^+_T(F),\ \Rot^-_T(F)$ and $\Rot_T(F),$ respectively.
\end{defi}

\subsection{$\sigma$-like maps}\label{subsec:sigmalikemaps}

Let $T\in\InfG$ and $F\colon T\to T$ a continuous map of degree $1$.
Define 
$$
\TR=\overline{\bigcup_{n\geq 0}F^n(\IR)}
$$
and
$X_F=\overline{T\setminus \TR}\cap r^{-1}([0,1))$.
Then $\TR\in \InfG$ (Lemma 5.2 in \cite{AR}), $X_F$ is
composed of finitely many finite graphs and $T=\TR\cup(X_F+\IZ)$.

\medskip
If $T$ is the lifting of the graph $\sigma$ (see Figure~\ref{fig:sigma}), 
then $X_F$ is, either empty, or an interval with an endpoint in $\TR$.
Maps with the same property will be called \emph{$\sigma$-like maps}.

\begin{defi}
Let $T\in\InfG$ and $F\colon T\to T$ a continuous map of degree $1$.
If $X_F$ is, either empty, or a nonempty interval such that 
$X_F\cap \TR$ is one of the endpoints of $X_F$,
we say that $F$ is a \emph{$\sigma$-like map} and we write $F\in\sigmap$.
\end{defi}

\begin{rem}
If $F$ is a $\sigma$-like map then so is $F^n$, because
$X_{F^n}\subset X_F$.
\end{rem}

This paper is devoted to the study of the rotation set of $\sigma$-like maps
when $X_F\not=\emptyset$. 
The study of the rotation set $\Rot_{\TR}(F)$ has already been 
done in \cite{AR}.

\section{Positive covering}
Let $F\in\sigmap$.
The interval $X_F$, when it is not empty, may be endowed with two
opposite orders. We choose the one such that $\min X_F$ is the
one-point intersection $X_F\cap \TR$. 
The retraction map $r_X\colon T\to X_F$ can be defined in a natural way 
by $r_X(x)=x$ if $x\in X_F$ and $r_X(x)=\min X_F$ if $x\in \TR$.

The notion of positive covering for subintervals of $\IR$
has been introduced in \cite{AR}. It can be extended for subintervals of
any subset of $T$ on which a retraction can be defined. In this
paper, we shall use positive covering on $X_F$.
All properties of positive covering remain valid in this context.
In particular, if a compact interval $I$ positively $F$-covers itself, 
then $F$ has a fixed
point in $I$ (Proposition~\ref{prop:+cover-periodic}).

\begin{defi}
Let $T\in\InfG$, $F\in\sigmap$, $I,J$ two non empty compact
subintervals of $X_F$, $n$ a positive integer and $p\in\IZ$. 
We say that $I$ {\em positively $F^n$-covers} $J+p$
and we write $I\pluscover{F^n} J+p$ if there exist $x\leq y$ in $I$ such that
$r_X(F^n(x)-p)\leq \min J$ and $\max J\leq r_X(F^n(y)-p)$. In this
situation, 
we also say that $I+q$ positively $F^n$-covers $J+p+q$ for all $q\in \IZ$.
\end{defi}

\begin{rem}
If $F^n(x)\in \TR$ and $J\subset X_F$, then the inequality 
$r_X(F^n(x)-p)\leq \min J$ is automatically satisfied. 
We shall often use this remark to prove that
an interval positively covers another.
\end{rem}

We introduce some definitions in order to handle sequences of
positive coverings.

\begin{defi}
Let $T\in\InfG$ and $F\in\sigmap$. If we have the following sequence of
positive coverings:
$$
\CC\colon I_0+p_0\pluscover{F^{n_1}}I_1+p_1\pluscover{F^{n_2}}I_2+p_2
\cdots\qquad \cdots
I_{k-1}+p_{k-1}\pluscover{F^{n_k}} I_k+p_k,
$$
(where $I_0,\ldots,I_k$ are non empty compact
subintervals of $X_F$, $n_1,\ldots,n_k$ are positive integers and $p_0,\ldots, 
p_k\in\IZ$), then $\CC$ 
is called a \emph{chain} of intervals for $F$. Its \emph{length} is
$L_F(\CC)=n_1+\cdots+n_k$, and its \emph{weight} is 
$W_F(\CC)=p_k-p_0$. A point $x$ \emph{follows} the chain $\CC$ if
$F^{n_1+\cdots+n_i}(x)\in I_i+p_i$ for all $0\leq i\leq k$.

If $i\in\IZ$, the chain $\CC+i$ is
the translation of $\mathcal C$, that is
$$
\CC+i\colon I_0+p_0+i\pluscover{F^{n_1}}I_1+p_1+i\pluscover{F^{n_2}}\cdots
\qquad\cdots I_{k-1}+p_{k-1}+i\pluscover{F^{n_k}} I_k+p_k+i.
$$

If $\CC'$ is another chain of intervals beginning with $I_k+p$ for some
$p\in\IZ$, then $\CC\CC'$ is the concatenation of 
$\CC$ and $(\CC'-p+p_k)$. 
If $I_k=I_0$, then $\CC^n$ is the $n$-times concatenation
$\CC\cdots \CC$ if $n\geq 1$ and $\CC^0$ is the empty chain.
\end{defi}

The next properties are straightforward.

\begin{lem}\label{lem:chain}
Let $T\in\InfG$ and $F\in\sigmap$.
\begin{itemize}
\item
If $\CC$ is a chain of intervals for $F^n$, then it is also
a chain of intervals for $F$ and $L_F(\CC)=nL_{F^n}(\CC)$ and $W_F(\CC)=
W_{F^n}(\CC)$. 
Since the weight is independent of the power of the map, we shall
denote it by $W(\CC)$.
\item
If $\CC,\CC'$ are two chains of intervals for $F$ that can be concatenated,
then $L_F(\CC\CC')=L_F(\CC)+L_F(\CC')$ and $W(\CC\CC')=W(\CC)+W(\CC')$.
\end{itemize}
\end{lem}

The next proposition is \cite[Proposition 2.3]{AR} (rewritten in some less 
general form). 

\begin{prop}\label{prop:+cover-periodic}
Let $T\in\InfG$, $F\in \sigmap$ and $\CC$ a chain of subintervals of $X_F$
such that $\CC$ starts with some interval $I_0$ and ends with a translation of
$I_0$ (i.e., $I_0+p$ for some $p\in\IZ$). Then there exists a point 
$x_0$ following the chain $\CC$ such that $F^{L_F(\CC)}(x_0)=x_0+W(\CC)$.
\end{prop}

The next lemma says that if two intervals $I,J$ both positively cover
translations of $I$ and $J$, then every
rational number in the interval corresponding to this ``horseshoe''
can be obtained as a rotation number of a periodic {\modi} point.
This will be a key tool.

\begin{lem}\label{lem:+horseshoe}
Let $T\in\InfG$, $G\in\sigmap$,
$I,J$ two non empty compact subintervals of $X_G$
and $m_1, m_2\in\IZ$ such that
\begin{eqnarray*}
I\pluscover{G}I+m_1&\mbox{ and }&I\pluscover{G}J+m_1,\\
J\pluscover{G}I+m_2&\mbox{ and }&J\pluscover{G}J+m_2.
\end{eqnarray*}
Suppose that $m_1\leq m_2$. For every $p/q\in [m_1,m_2]$,
there exists $\CC$ a chain of intervals for $G$ in which all the intervals
are translations of $I$ and $J$, 
and $p/q=W(\CC)/L_G(\CC)$. Moreover, there exists a periodic
{\modi} point $x\in I\cup J$ such that $\rhos[G](x)=p/q$.

If $p/q\not=m_2$, then $\CC$ can be chosen such that
the first interval is $I$ and the last interval is a translation of $I$,
and the periodic {\modi} point $x$ can be chosen in $I$.
\end{lem}

\begin{demo}
By considering $G-m_1$ instead of $G$, we may suppose that $m_1=0$
(use Lemma~\ref{lem:first-properties}).
Since $p/q\in [0,m_2]$, we have $0\leq p\leq m_2q$. 
If $p/q=m_2$, we take $\CC\colon J\pluscover{G}J+m_2$. By 
Proposition~\ref{prop:+cover-periodic}, there exists a point $x\in J$ such
that $F(x)=x+m_2$, and hence $\rho_G(x)=m_2$.

If $p=0$, we take $\CC\colon I\pluscover{G} I$.
If $1\leq p\leq m_2q-1$, we take 
$$
\CC\colon (I\pluscover{G}I)^{m_2q-1-p}
(I\pluscover{G}J)(J\pluscover{G}J+m_2)^{p-1}(J\pluscover{G}I+m_2)
$$
In these two cases, it is straightforward that $W(\CC)/L_G(\CC)=p/q$.
By Proposition~\ref{prop:+cover-periodic}, there exists a point $x\in I$
such that $G^{L_G(\CC)}(x)=x+W(\CC)$, and so $x$ is periodic {\modi}
and $\rhos[G](x)=W(\CC)/L_G(\CC)=p/q$.
\end{demo}

\section{Study of the rotation set of $F$}

Let $T\in\InfG$ and $F\in\sigmap$.
Since $T=\TR\cup (X_F+\IZ)$, it is clear that
$\Rot(F)=\Rot_{\TR}(F)\cup \Rot_{X_F}(F)$, and the same holds with
$\Rot^+$ and $\Rot^-$. The rotation set $\Rot_{\TR}(F)$ has been studied
in \cite{AR}. 
Consequently, it remains to study the rotation set $\Rot_{X_F}(F)$.
The next theorem summarises the properties of $\RotR(F)$
(see Theorems 3.1, 3.11, 5.7 and 5.18 in \cite{AR}).

\begin{theo}\label{theo:rotR}
Let $T\in \InfG$ and $F\colon T\to T$ a continuous map of degree $1$. Then
$\Rot_{\IR}(F)$ is a non empty compact interval and, if $\TR$ is defined
as above,
$\Rot_{\TR}(F)=
\Rot^+_{\TR}(F)=\Rot^-_{\TR}(F)=\Rot_{\IR}(F)$. Moreover, if $r\in 
\Rot_{\IR}(F)\cap\IQ$, then there exists a periodic {\modi} point
$x\in \TR$ such that $\rhos(x)=r$.
\end{theo}

\subsection{Partition of $X_F$}

If $F^n(x)\in \TR$ for some $n$, then $\rhos(x)\in\Rot_{\TR}(F)$.
We already know the properties of $\Rot_{\TR}(F)$.
Therefore, it is sufficient to consider the points $x\in X_F$ whose
orbit does not fall in $\TR$, or equivalently the points
in $X_{\infty}=\{x\in X_F\mid \forall n\geq 1, F^n(x)\in X_F+\IZ\}$.

Our first step consists in dividing $X_F$ according to the translations
of the images with respect to $X_F+\IZ$. If $F(x)\in X_F+p$
and $F(y)\in X_F+p'$ with $p\not=p'$, then necessarily there is a 
gap between $x$ and $y$ by continuity. Thus we can include the points
$\{x\in X_F\mid F(x)\in X_F+\IZ\}$ in a finite union of disjoint
compact intervals such that, for each $I$ among these intervals, there
is a unique integer $p$ satisfying $F(I)\cap (X_F+p)\not=\emptyset$.

\begin{lem}\label{lem:Xi}
Let $T\in\InfG$ and $F\in\sigmap$. There exist an integer $N\geq 0$, 
non empty compact subintervals $X_1,\ldots, X_N$ of $X_F$ and
integers $p_1,\ldots, p_N$ in $\IZ$ such that
\begin{enumerate}
\item $X_1<X_2<\cdots<X_N$,
\item $F(X_i)\subset (X_F+p_i)\cup \TR$ for all $1\leq i\leq N$,
\item $F(\min X_i)=\min X_F+p_i$ for all $1\leq i\leq N$,
\item $p_{i+1}\not=p_i$ for all $1\leq i\leq N-1$,
\item $F(X_F\setminus (X_1\cup\cdots\cup X_N))\cap (X_F+\IZ)=\emptyset$.
\end{enumerate}
\end{lem}

\begin{demo}
If $F(X_F)\cap(X_F+\IZ)=\emptyset$, we take $N=0$ and there is nothing to do.
Otherwise, we can define $a_1=\min\{x\in X_F\mid F(x)\in X_F+\IZ\}$
and $p_1\in\IZ$ such that $F(a_1)\in X_F+p_1$. We define
$$
b_1=\max\{x\in [a_1,\max X_F]\mid F(x)\in X_F+p_1\mbox{ and }
F([a_1,x])\subset (X_F+p_1)\cup \TR\},
$$
and $X_1=[a_1,b_1]$. Then $X_1$ satisfies (ii). Moreover, $F(\min X_F)\in\TR$
because $\min X_F\in\TR$, which implies that $F([\min X_F,a_1])$ contains
$\min X_F+p_1$. Thus $F(a_1)=\min X_F+p_1$ by minimality of $a_1$, which is
(iii) for $X_1$.

We define $X_2,\ldots, X_M$ inductively.
Suppose that $X_i=[a_i,b_i]$ and $p_i$ are already defined and that
$b_i$ verifies: 
$$
b_i=\max\{x\in [a_i,\max X_F]\mid F(x)\in X_F+p_i\mbox{ and }
F([a_i,x])\subset (X_F+p_i)\cup \TR\}.
$$
If $F((b_i,\max X_F])\cap(X_F+\IZ)=\emptyset$, then we take $N=i$ and the construction is
over. Otherwise, we define
\begin{equation}\label{eq:ai}
a_{i+1}=\inf\{x\in (b_i,\max X_F]\mid F(x)\in X_F+\IZ\}.
\end{equation}
We first show that $a_{i+1}$ is actually defined by a minimum in~\eqref{eq:ai}.
By definition, there exists a sequence of points $x_n\in (b_i,\max X_F]$
tending to $a_{i+1}$ and such that $F(x_n)\in X_F+\IZ$. Let $m_n\in\IZ$
such that $F(x_n)\in X_F+m_n$.
By continuity, $\lim_{n\to+\infty} r\circ F(x_n)=r\circ F(a_{i+1})$.
Since $r\circ F(x_n)=r(\min X_F)+m_n$, 
this implies that the sequence of integers $(m_n)_{n\geq 0}$ is
ultimately constant, and equal to some integer $p_{i+1}$. Then $F(a_{i+1})
=\lim_{n\to+\infty}F(x_n)\in X_F+p_{i+1}$. 
By continuity, $F([a_{i+1}, x_n])\subset (X_F+p_{i+1})\cup \TR$ for 
all $n$ large enough.
Moreover, $F((b_i,a_{i+1}))\cap(X_F+\IZ)=\emptyset$ by definition of $a_{i+1}$.
If $p_{i+1}=p_i$, then, for $n$ large enough, we would have 
$$
F(x_n)\in X_F+p_i\mbox{ and }
F([b_i,x_n])\in (X+p_i)\cup \TR,
$$
which would contradict the definition of $b_i$ because $x_n>b_i$. 
Hence $p_{i+1}\not=p_i$.
This implies that $a_{i+1}>b_i$.
Since $F((b_i,a_{i+1}))$ is non empty and included in $\TR$, necessarily 
$F(a_{i+1})$ is equal to $\min X_F+p_{i+1}$ by minimality of $a_{i+1}$.

Finally, we define
$$
b_{i+1}=\max\{x\in [a_{i+1},\max X_F]\mid F(x)\in X_F+p_{i+1}\mbox{ and }
F([a_{i+1},x])\subset (X_F+p_{i+1})\cup \TR\},
$$
and $X_{i+1}=[a_{i+1},b_{i+1}]$. Then
$X_{i+1}>X_i$ and (ii), (iii) and (iv) are satisfied.

By uniform continuity of $r\circ F$ on the compact set $X_F$, 
there exists $\delta>0$ such that, if $x,y\in X_F$
with $|x-y|<\delta$, then $|r\circ F(x)-r\circ F(y)|<1$. This
implies that $|a_{i+1}-b_i|\geq\delta$, which ensures that 
the number of intervals $X_i$ is finite, and the construction
ultimately stops. By construction, (v) is satisfied.
\end{demo}

\begin{rem}

\begin{enumerate}
\item
The fact that the sets $X_1,\ldots, X_N$ are intervals
is very important because it will allow us to use positive coverings.
Note that we cannot ask that $F(X_i)\subset (X_F+p_i)$, even if we do not
require that $p_{i+1}\neq p_i$. Indeed, if $\min X_F$ is a fixed point,
the map $F$ may oscillate infinitely many times between $X_F$ and $\TR$
in any neighbourhood of $\min X_F$, and in this case the number of 
connected components of $F(X_F)\cap X_F$ is infinite.
\item
In the partition of $X_F$ into $X_1,\ldots,X_N,X_F\setminus(X_1\cup\cdots
\cup X_N)$, the set $X_F\setminus(X_1\cup\cdots \cup X_N)$ plays the role of
``dustbin'', and we can code the itinerary of every point in $X_{\infty}$
with respect to $X_1,\ldots,X_N$. More precisely, if $F^n(x)\in
X_F\setminus(X_1\cup\cdots \cup X_N)+\IZ$, then $x\not\in X_{\infty}$.
Therefore, for every $x\in X_{\infty}$, $\forall n\geq 0,\exists ! \omega_n
\in\{1,\ldots,N\}$ such that $F^n(x)\in X_{\omega_n}+\IZ$. The rotation
number of $x$ can be deduced from this coding sequence because
$\forall n\geq 0$, $F^n(x)\in X_{\omega_n}+p_{\omega_0}+\cdots+p_{\omega_{n-1}}$.
\item
It can additionally been shown that
$F(\max X_i)=\min X_F+p_i$ for all $1\leq i\leq N-1$ and, for $i=N$,
either $F(\max X_N)=\min X_F+p_N$, or $\max X_N=\max X_F$.
\end{enumerate}
\end{rem}

\subsection{Periodic {\modi} points associated to the endpoints of
rotation sets}

When proving that every rational number in the rotation set is the rotation
number of a periodic {\modi} point, we shall make a distinction
between rational numbers in the interior and rational numbers in the boundary
(the same distinction is necessary to deal with $\Rot_{\IR}(F)$ \cite{AR}).
For rational numbers in the boundary, harder to handle, we shall need the
following lemma, which is the analogous of \cite[Lemma~5.15]{AR} in our
context. This lemma is aimed to be applied first with $T'=\TR$, $Y=X_F$
and $Z=X_1$, where $X_1$ is defined in
Lemma~\ref{lem:Xi}. After dealing
with $X_1$, an induction will be done to deal with $X_2,\ldots, X_N$, that
is why the lemma is stated with general notations.

\begin{lem}\label{lem:minRot>0}
Let $T\in\InfG$ and $F\in\sigmap$. Let $T'$ be a closed connected 
subset of $T$ such that $\TR\subset T'$, $T'+1=T'$ and $F(T')\subset 
T'$.
Let $Y$ denote the compact subinterval of $X_F$ equal to 
$\overline{T\setminus T'}\cap r^{-1}([0,1))$ and define
$Y_{\infty}=\{x\in Y\mid \forall n\geq 1, F^n(x)\in Y+\IZ\}$.
Let $Z$ be a compact subinterval of $Y$ such that $F(\min Z)=\min Y$ and 
$F(Z)\cap Z\neq\emptyset$. Assume that
$\inf \Rot_{Z\cap Y_{\infty}}(F)\geq 0$ and
$$
\forall x\in \overline{\bigcup_{n\geq 0} (F^n(Z)+\IZ)}\cap Y_{\infty},
\forall n\geq 1, F^n(x)\not=x.
$$
Then $\inf \Rot_{Z\cap Y_{\infty}}(F)> 0$.
\end{lem}

\begin{demo}
We note $x_0=\min Y$. By assumption, $F^n(x_0)\in T'$ for all $n\geq 0$.
Let $$Y'=\overline{\bigcup_{n\geq 0} (F^n(Z)+\IZ)}\cap Y.$$
For all $n\geq 1$, $F^n(x_0)\in T'$, which implies
that $(F^n(Z)+\IZ)\cap Y$ is an interval containing $x_0$.
In addition, $F(Z)\cap Z\neq\emptyset$ by assumption. This implies
that $Y'$ is a compact subinterval of $Y$ containing $x_0$.

Let $x\in Y'$ and $k\in\IZ$ such that $F(x)\in Y+k$.
Suppose that $F(x)-k\geq x$. Then $[x_0,x]\pluscover{F}[x_0,x]+k$ and,
according to Proposition~\ref{prop:+cover-periodic}, there exists
$y\in [x_0,x]\subset Y'$ such that $F(y)=y+k$. Hence $y\in Y_{\infty}$
and $\rhos(y)=k$. By assumption, this is possible  only if $k>0$.
Hence
\begin{equation}\label{eq:k<0}
\mbox{if }x\in Y' \mbox{ and }F(x)\in Y+k \mbox{ with } 
k\leq 0, \mbox{ then } F(x)-k<x.
\end{equation}
We are going to prove that there exists an integer $N_0$ such that
\begin{equation}\label{eq:N0}
\mbox{if }y\in Y' \mbox{ verifies }
\forall\; 1\leq n\leq N, \exists\, k_n\leq 0,  F^n(y)\in Y+k_n,
\mbox{ then }N\leq N_0.
\end{equation}
Let $y$ verify the assumption of (\ref{eq:N0}) and note $k_0=0$.
By Equation~(\ref{eq:k<0}), $(F^i(y)-k_i)_{0\leq i\leq N}$ is a decreasing
sequence in $Y$. Since $Y'$ is an interval containing $y$ and $\min Y$,
all the points of this sequence belong to $Y'$ too.

According to \eqref{eq:k<0},
for all $k\leq 0$ and all $x\in Y'$, $F(x)-k\not=x$. If $d$ denotes a
distance on $T$, then
$$
\forall k\leq 0,\ \delta_k=\inf\{d(x,F(x)-k)\mid x\in Y'\}>0
$$
because $Y'$ is compact. Moreover, the set of integers $k$ such that
$F(Y')\cap (Y+k)\not=\emptyset$ is finite, and so
$$
\delta=\inf_{k\leq 0} \delta_k>0.
$$ 
Consequently, 
$(F^i(y)-k_i)_{0\leq i\leq N}$ is a decreasing sequence in $Y'$ and,
for all $0\leq i\leq N-1$,
$d(F^{i+1}(y)-k_{i+1}, F^i(y)-k_i)\geq \delta$.
This implies that $\diam{Y'}\geq N\delta$. This proves that
(\ref{eq:N0}) holds if $N_0\geq\frac{\diam{Y'}}{\delta}$.

Now, let $x\in Z\cap Y_{\infty}$. For all $n\geq 0$, $F^n(x)\in Y'+\IZ$.
According to (\ref{eq:N0}), there exists an increasing sequence of
positive integers $(n_i)_{i\geq 1}$ and integers $(k_i)_{i\geq 1}$
such that 
$$
\forall i\geq 1,\ F^{n_i}(x)\in Y'+k_i,\
n_{i+1}-n_i\leq N_0 \quad\mbox{and}\quad k_{i+1}\geq k_i+1.
$$
This implies that
$\orhos(x)\geq 1/N_0>0$. This concludes the proof of the lemma.
\end{demo}

Let us restate Lemma~\ref{lem:minRot>0} when $T'=\TR$, $Y=X_F$ and 
$Z=X_1$: if $\inf \Rot_{X_1\cap X_{\infty}}(F)=0$ and
$F(X_1)\cap X_1\neq\emptyset$,  then there exist
a point $x\in \overline{\bigcup_{n\geq 0} (F^n(X_1)+\IZ)}\cap X_{\infty}$ 
and an integer $n\geq 1$ such that $F^n(x)=x$. 
Note that when $p_1=0$, the assumption $F(X_1)\cap X_1\neq\emptyset$ 
is fulfilled as soon as
$\Rot_{X_1\cap X_{\infty}}(F)\neq\emptyset$. Indeed, if $F(X_1)\cap X_1=
\emptyset$, then $F(X_1)\cap (X_F+\IZ)\subset [\min X_F,\min X_1)$,
and thus $X_1\cap X_{\infty}=\emptyset$.
More generally, this result implies that, if 
$\inf \Rot_{X_1\cap X_{\infty}}(F)=p/q$, then there exists a 
periodic {\modi} point $x\in X_F$ such that $\rho_F(x)=p/q$
(just apply what precedes to the map $F^q-p$). 
This is the analogous in our context of \cite[Theorem~5.18]{AR}.

\subsection{Rotation set of $X_1$}

In the sequel, we shall heavily use the fact that $X_F$ is an interval
with an endpoint in $\TR$. 
By definition, $\min X_F$ belongs to $\TR$ and $\TR$ is invariant
by $F$. Hence $F(\min X_F)\in \TR$. Therefore, if $I$ is a subinterval of
$X_F$ such that
$\min I=\min X_F$
and $F(I)\cap X_F\not=\emptyset$, then necessarily $F(I)\cap X_F$
is an interval containing $\min X_F$. This simple observation allows us
to study the rotation set of the interval $X_1$ defined in
Lemma~\ref{lem:Xi}. This is done in Proposition~\ref{prop:key}, 
by considering $T'=\TR$ and $Y_1,\ldots, Y_M=X_1,\ldots, X_N$.
When this is done for $X_1$, the idea is to proceed by induction for the
rotation sets of $X_2,\ldots, X_N$, this is why the proposition is stated with
more general notations.

In the proof of Proposition~\ref{prop:key}, 
we shall need the next, technical lemma.

\begin{lem}\label{lem:limsup}
Let $(n_k)_{k\geq 0}$ be a sequence of real numbers bounded from above 
by some constant~$C$.
Let $\eps>0$, $$L=\limsup_{k\to+\infty}\frac{n_0+\cdots+n_{k-1}}{k}$$
and $l<L-\eps$. 
Then there exists an integer $k\geq 1$
such that $\frac{n_0+\cdots+n_{k-1}}{k}\geq L-\eps$ and $n_k\not=l$.
\end{lem}

\begin{demo}
Let $K$ be an integer such that $C/K<\eps/2$. Let
$$
E=\left\{k\in\IZ^+\;\big|\; \frac{n_0+\cdots+n_k}{k+1}\geq L-\eps/2\right\}.
$$
The set $E$ is infinite by definition of $L$. 
We are going to do a proof by absurd. We assume that 
\begin{equation}\label{eq:nk=l}
n_k=l\ \mbox{ for all $k\in E$ such that }k>K.
\end{equation}
If $E$ contains all integers $n\geq N$ for some $N$, then $L=l$,
which is absurd. Thus there exists an integer $k\geq K$ such that $k\in E$ and
$k-1\not\in E$. We have
$$
\frac{n_0+\cdots+n_k}{k+1}=\frac{k}{k+1}\cdot\frac{n_0+\cdots+n_{k-1}}{k}+
\frac{1}{k+1} n_k.
$$
By definition of $E$, $\frac{n_0+\cdots+n_k}{k+1}\geq L-\eps/2$
and $\frac{n_0+\cdots+n_{k-1}}{k}<L-\eps/2$. Moreover, $n_k=l<L-\eps/2$
by \eqref{eq:nk=l}. Thus
$$
\frac{n_0+\cdots+n_k}{k+1}< \frac{k}{k+1} (L-\eps/2) +\frac{1}{k+1}(L-\eps/2)
=L-\eps/2,
$$
which is a contradiction.
Therefore, \eqref{eq:nk=l} does not hold, and
there exists $k\in E$ such that $k> K$ and $n_k\not=l$.
Moreover, 
\begin{eqnarray*}
\frac{n_0+\cdots+n_{k-1}}{k}&=&
\frac{k+1}{k}\cdot \frac{n_0+\cdots+n_k}{k+1}-\frac{n_k}{k}\\
&>& \frac{n_0+\cdots+n_k}{k+1} -\frac{C}{K}\\
&> & L-\eps/2 -\eps/2 =L-\eps
\end{eqnarray*}
Such an integer $k$ is suitable.
\end{demo}

\begin{prop}\label{prop:key}
Let $T\in\InfG$ and $F\in\sigmap$. Let $T'$ be a connected 
subset of $T$ such that $\TR\subset T'$, $T'+1=T'$ and $F(T')\subset T'$.
Let $Y$ denote the compact subinterval of $X_F$ equal to 
$\overline{T\setminus T'}\cap r^{-1}([0,1))$ and define
$Y_{\infty}=\{x\in Y\mid \forall n\geq 1, F^n(x)\in Y+\IZ\}$.
Let $Y_1,\ldots, Y_M$ be  compact subintervals of $Y$ and
$q_1,\ldots, q_M\in\IZ$ such that:
\begin{itemize}
\item[a)] $Y_1<\cdots<Y_M$,
\item[b)] $F(Y_i)\subset (Y+q_i)\cup T'$ for all $1\leq i\leq M$,
\item[c)] $F(\min Y_1)=\min Y\ {\modi}$,
\item[d)] $F(Y\setminus (Y_1\cup\cdots\cup Y_M))\cap (Y+\IZ)=\emptyset$.
\end{itemize}
Assume that $Y_1\cap Y_{\infty}\not=\emptyset$. Then there exists a compact
interval $I\subset \IR$ such that:
\begin{enumerate}
\item
$\Rot_{Y_1\cap Y_{\infty}}(F)=\Rot^+_{Y_1\cap Y_{\infty}}(F)=
\Rot^-_{Y_1\cap Y_{\infty}}(F)=I$,
\item
$I$ contains $q_1$,
\item
there exists $a\in Y_1\cap Y_{\infty}$ such that $F(a)=a+q_1$ and $[\min Y_1,a)\cap Y_{\infty}
=\emptyset$,
\item
if $r\in \Int{I}\cap\IQ$ then there exists a periodic {\modi} point
$x\in Y_1\cap Y_{\infty}$ with $\rhos(x)=r$.
\item if $r\in\partial I\cap\IQ$ then there exists a periodic {\modi}
point $x\in \overline{\bigcup_{n\geq 0} (F^n(Y_1)+\IZ)}\cap Y_{\infty}$ 
with $\rhos(x)=r$.
\end{enumerate}
\end{prop}

\begin{demo}
We first prove (iii) under an additional assumption:
\begin{equation}\label{eq:iii}
\mbox{If }\exists y\in Y_1\mbox{ such that }\forall n\geq 1, F^n(y)\in Y_1+\IZ,
\mbox{ then (iii) holds.}
\end{equation}

Let $G=F-q_1$ and $a_0=\min Y_1$. 
Then $G(Y_1)\subset Y\cup T'$. For all $n\geq 0$, 
$G^n(y)\in Y_1$, and in particular $G^n(y)\geq a_0$. 
We define inductively  a sequence of points $(a_i)_{i\geq 1}$ such
that $a_i\in [a_{i-1},y]$, $G^i(a_i)=a_0$ and 
$[a_0,a_i)\cap Y_\infty=\emptyset$ for all $i\geq 1$.
\begin{itemize}
\item Since $G(a_0)=\min Y$ (by assumption (c)) and $G(y)\geq a_0$, we have 
$a_0\in G([a_0,y])$ by continuity. Thus there exists $a_1\in [a_0,y]$ such that
$G(a_1)=a_0$. We choose $a_1$ minimum with this property, which implies
that $G([a_0,a_1))\cap Y_1=\emptyset$. Hence $[a_0,a_1)\cap Y_\infty=
\emptyset$.
\item Assume that $a_0,\ldots, a_i$ are already defined. 
Since $G^{i+1}(a_i)=G(a_0)=\min Y$ and $G^{i+1}(y)\geq a_0$, the point $a_0$
belongs to $G^{i+1}([a_i,y])$ by continuity. Thus there exists $a_{i+1}
\in [a_i,y]$ such that $G^{i+1}(a_{i+1})=a_0$. 
We choose $a_{i+1}$ minimum with this property, which implies
that $G^{i+1}([a_i,a_{i+1}))\cap Y_1=\emptyset$. Hence
$[a_i,a_{i+1})\cap Y_\infty=\emptyset$. This concludes
the construction of $a_{i+1}$.
\end{itemize}
The sequence $(a_i)_{i\geq 0}$ is non decreasing and contained in the
compact interval $Y_1$. Therefore $a=\lim_{i\to+\infty}a_i$ exists
and belongs to $Y_1$. Since $G(a_{i+1})=a_i$, we get that $G(a)=
a$. In other words, $F(a)=a+q_1$. This implies
that $a\in Y_{\infty}$.  Moreover, $[a_0,a)=\bigcup_{i\geq 0}
[a_i,a_{i+1})$, and thus
$[a_0,a)\cap Y_\infty=\emptyset$. This proves \eqref{eq:iii}.

\medskip
We split the rest of the proof into two cases.

\textbf{Case 1:} $F(Y_1)\cap(Y_i+\IZ)=\emptyset$ for all $i\geq 2$ (this 
includes the case $M=1$).\\
Then $F(Y_1\cap Y_{\infty})\subset Y_1+q_1$ and 
$F^n(Y_1\cap Y_{\infty})\subset Y_1+nq_1$. Thus, for all 
$x\in Y_1\cap Y_{\infty}$, the rotation number $\rhos(x)$ exists and is equal
to $q_1$. We take $I=\{q_1\}$ and we get (i) and (ii).

Since $Y_1\cap Y_{\infty}$ is not empty, there exists
a point $y$ such that $\forall n\geq 0$, $F^n(y)\in Y_1+\IZ$.
Then \eqref{eq:iii} gives (iii), which implies
(v) in the present case, and (iv) is empty. 

\medskip
\textbf{Case 2:} there exists $i\geq 2$ such that $F(Y_1)\cap(Y_i+\IZ)\not=
\emptyset$.\\
Since $Y_1<\cdots<Y_M$ and $F(Y_1)\subset (Y+q_1)\cup T'$, there
exists $x\in Y_1$ such that $F(x)\in Y+q_1$ and $F(x)\geq \max (Y_1+q_1)$.
Moreover, $F(\min Y_1)=\min Y+q_1$ by assumption. This implies that
\begin{equation}\label{eq:Z1+cover}
Y_1\pluscover{F} Y_1+q_1.
\end{equation}
By Proposition~\ref{prop:+cover-periodic}, there exists $y\in Y_1$ 
such that $F(y)=y+q_1$. Thus we can use \eqref{eq:iii} to get (iii).

Let $x\in Y_{\infty}$. By assumption d), for all $n\geq 0$, there exists 
$\omega_n\in\{1,\ldots,M\}$ such that $F^n(x)\in Y_{\omega_n}+\IZ$.
The sequence $(\omega_n)_{n\geq 0}$ is called the \emph{itinerary}
of $x$. The next two results are straightforward.
\begin{equation}
\forall n\geq 0,\ F^n(x)\in Y_{\omega_n}+q_{\omega_0}+q_{\omega_1}+\cdots+q_{\omega_{n-1}}. \label{eq:itinerary-Fn}
\end{equation}
\begin{eqnarray}
&\displaystyle \orhos(x)=\limsup_{n\to+\infty} \frac{q_{\omega_0}+q_{\omega_1}+\cdots+q_{\omega_{n-1}}}{n}\quad\mbox{and}\quad
\urhos(x)=\liminf_{n\to+\infty} \frac{q_{\omega_0}+q_{\omega_1}+\cdots+q_{\omega_{n-1}}}{n};\label{eq:itinerary-rho}&\\
&\mbox{if the limit exists, it is }\rhos(x).&\nonumber
\end{eqnarray}

Let $S=\sup\Rot^+_{Y_1\cap Y_{\infty}}(F)$. 
Necessarily, $S\geq q_1$ because $q_1\in \Rot_{Y_1\cap Y_\infty}(F)$ by (iii).
We are going to show
that $[q_1,S]\subset \Rot_{Y_1\cap Y_{\infty}}(F)$.
If $S=q_1$ there is nothing to prove, and so we suppose that $S>q_1$.
Let $k$ be an integer such that $S>q_1+1/k$. 
Let $y_k$ be a point in $Y_1\cap Y_{\infty}$
such that $\orhos(y_k)\geq S-1/2k$, and let $(\omega_n)_{n\geq 0}$ denote
the itinerary of $y_k$. By (\ref{eq:itinerary-rho}),
$\limsup_{n\to+\infty}\frac{q_{\omega_0}+q_{\omega_1}+\cdots+q_{\omega_{n-1}}}{n}=\orhos(y_k)$. Applying Lemma~\ref{lem:limsup} with $L=\orhos(y_k)$, 
$l=q_1$ and $\eps=1/2k$, we get that there exists an integer $n$ such that
\begin{equation}\label{eq:NK/nk}
\frac{q_{\omega_0}+q_{\omega_1}+\cdots+q_{\omega_{n-1}}}{n}\geq \orhos(y_k)-
1/2k\geq S-1/k>q_1 \mbox{ and } \omega_n\not=1.
\end{equation}
By (\ref{eq:itinerary-Fn}), $F^n(y_k)\in Y_{\omega_n}+q_{\omega_0}+
q_{\omega_1}+\cdots+q_{\omega_{n-1}}$. Since $\omega_n\not=1$, we have 
$Y_{\omega_n}>Y_1$, and thus $F^n(y_k)-(q_{\omega_0}+
q_{\omega_1}+\cdots+q_{\omega_{n-1}})>\max Y_1$.
Moreover, $F^n(\min Y_1)\in F^{n-1}(\overline{T'})\subset \overline{T'}$ 
by assumption c) and
invariance of $T'$. If we let $I_k=[\min Y_1, y_k]$, 
$N_k=q_{\omega_0}+\cdots+q_{\omega_{n-1}}$ and $n_k=n$, we have then:
\begin{equation}\label{eq:Ik+cover}
I_k\pluscover{F^{n_k}}Y_1+N_k.
\end{equation}
Since $I_k\subset Y_1$, Equations (\ref{eq:Z1+cover}) and 
(\ref{eq:Ik+cover}) give:
\begin{equation}\label{eq:Ik+Y1}
\begin{array}{rcl}
I_k\pluscover{F^{n_k}}I_k+N_k&\mbox{ and }&I_k\pluscover{F^{n_k}}Y_1+N_k,\\
Y_1\pluscover{F^{n_k}}I_k+n_kq_1&\mbox{ and }&Y_1\pluscover{F^{n_k}}Y_1+n_kq_1.
\end{array}
\end{equation}

Let $r\in [q_1,S)\cap\IQ$. By (\ref{eq:NK/nk}), 
there exists an integer $k$ such that
$\frac{N_k}{n_k}> r$. We apply Lemma~\ref{lem:+horseshoe} with 
$I=Y_1$, $J=I_k$, $G=F^{n_k}$, $m_1=n_kq_1$, $m_2=N_k$ and 
$p/q=rn_k\in[q_1n_k,N_k)$: 
\begin{gather}\label{eq:Ci}
\mbox{$\exists\;\CC_r$ a chain of intervals for
$F$ whose first and last intervals are translations of $Y_1$,} \nonumber \\
\mbox{such that }r=W(\CC_r)/L_F(\CC_r),
\end{gather} 
and 
\begin{equation}\label{eq:rhox=r}
\mbox{there exists a periodic {\modi} point
$x\in Y_1$ with }\rho_F(x)=\frac{1}{n_k}\rho_{F^{n_k}}(x)=r.
\end{equation} 
We need
to show that $x\in Y_\infty$. This is a consequence of the following fact.

\medskip\noindent\textbf{Fact.} 
Let $x\in Y_1$ such that $F^n(x)\in Y_1+\IZ$ for
infinitely many $n$. Then, either $x\in Y_{\infty}$, or there exists
$n$ such that $F^n(x)$ is a fixed {\modi} point in $Y_1\cap Y_{\infty}\cap T'$.

\medskip\noindent\textit{Proof of the fact.}
If $\forall n\geq 0,\ F^n(x)\not\in T'$, then $x\in Y_{\infty}$ by assumption
d) of the proposition. Suppose on the contrary that 
there exists $n_0$ such that $F^{n_0}(x)\in T'$. Hence
$F^n(x)\in T'$ for all $n\geq n_0$.
Let $e=\min Y$. By definition of $Y$, the set $Y\cap T'$ is included in 
$\{e\}$ (we have not supposed that $T'$ is closed, and thus $Y\cap T'$ may
be empty). Note that $Y_1\cap T'$ is empty if $\min Y_1>e$.
By assumption, there exists $n_1\geq n_0$ such that $F^{n_1}(x)\in Y_1+\IZ$. 
Hence $F^{n_1}(x)\in (Y_1+\IZ)\cap T'$. This implies that $F^{n_1}(x)$ is
equal to $e{\modi}$, $\min Y_1=e\in T'$, and $e\in T'$. 
By assumption c) of the proposition,
$F(\min Y_1)=e {\modi}$. Thus, $e$ 
is a fixed {\modi} point in $Y_1$, and so $e\in Y_{\infty}$. This
ends the proof of the fact. $\Box$

\medskip
Now, let $\alpha\in [q_1,S]$. To show that there
exists $x\in Y_1\cap Y_{\infty}$ with $\rhos(y)=\alpha$, 
we use the same method as in the proof of \cite[Theorem 3.7]{AR}.
We choose a sequence of rational numbers
$r_i$ in $[q_1,S)\cap\IQ$ such that $\lim_{i\to+\infty}r_i=\alpha$
for all $i\geq 1$. 
Let $\CC_{r_i}$ be the chain of intervals given by \eqref{eq:Ci}.
We define
$$
\CD_n=(\CC{r_1})^{i_1}(\CC_{r_2})^{i_2}\cdots(\CC_{r_n})^{i_n}.
$$
Let $A_n$ be the set of points that follow the chain $\CD_n$. This set
is compact by definition, and it is 
not empty because it contains at least a periodic {\modi} point by
Proposition~\ref{prop:+cover-periodic}. Moreover, $A_{n+1}\subset A_n$.
Therefore $A=\bigcap_{n\geq 1}A_n \neq \emptyset$. 
In the proof of \cite[Theorem 3.7]{AR}, it is shown that,
if the sequence $(i_n)_{n\geq 1}$ increases sufficiently fast
and $(|r_{i_n}-\alpha|)_{n\geq 1}$ is non decreasing
then, for all $x\in A$, $\rhos(x)=\alpha$.

Moreover, the fact above implies that, for every $x\in A$ there exists
$n$ such that $F^n(x)\in Y_1\cap Y_{\infty}$. Obviously,
$\rhos(x)=\rhos(F^n(x))$.
This proves that $[q_1,\sup \Rot^+_{Y_1\cap Y_{\infty}}(F)]$ is included in
$\Rot_{Y_1\cap Y_{\infty}}(F)$; in addition, (iv) holds for all rational numbers
$r\in[q_1,\sup \Rot^+_{Y_1\cap Y_{\infty}}(F))$ by \eqref{eq:rhox=r}.
We can apply the same method to $[\inf \Rot^-{Y_1\cap Y_{\infty}}(F),q_1]$.
Finally, if we define
$I=[\inf \Rot^-_{Y_1\cap Y_{\infty}}(F),
\sup\Rot^+_{Y_1\cap Y_{\infty}}(F)]$, we get that
$I\subset\Rot_{Y_1\cap Y_{\infty}}(F)$, $q_1\in I$ (which is (ii)) 
and (iv) holds
for all $r\in \Int{I}\cap\IQ$. Since $\Rot^+_{Y_1\cap Y_{\infty}}(F)$
and $\Rot^-_{Y_1\cap Y_{\infty}}(F)$ both 
contain $\Rot_{Y_1\cap Y_{\infty}}(F)$ and are included in $I$, this
gives (i). 

Now we prove (v) for $\min I$ (the case with the maximum is symmetric, and
$\partial I$ is reduced to two points).
Suppose that $\min \Rot_{Y_1\cap Y_{\infty}}(F)=p/q$ and let $G=F^q-p$.
Then $\min \Rot_{Y_1\cap Y_{\infty}}(G)=0$. Note that $x\in Y_{\infty}+\IZ
\Rightarrow \forall n\geq 0,\ G^n(x)\in Y+\IZ$.
We apply Lemma~\ref{lem:minRot>0} to the map $G$ with $Z=Y_1$
and $\overline{T'}$. Since (iii) is fulfilled, the set $G(Z)\cap Z$ is
not empty.
By refutation of Lemma~\ref{lem:minRot>0}, we get that there exist $x\in 
\overline{\bigcup_{n\geq 0}(G^n(Y_1)+\IZ)}\cap Y_{\infty}$ and $n\geq 1$
such that $G^n(x)=x$. Therefore, $x$ is periodic {\modi} for $F$
and $\rhos(x)=p/q$. Moreover, $x\in Y_\infty$ according to the fact above.
This gives (v) and concludes the proof of the proposition.
\end{demo}

\begin{ex}\label{ex:counterex}
The periodic {\modi} point $x$ given by Proposition~\ref{prop:key}(v)
may not be in $Y_1$. 

Let $T$ be the universal covering of the graph sigma, and let $F\colon 
T\to t$ be the continuous map of degree 1 such that $F|_{\IR}=Id$
and $F$ is defined on the branch of $T$ by:
\begin{itemize}
\item
$F(a)=F(c)=a$, $F(b)=e$, $F(d)=a+1$, $F(e)=e+1$,
\item $F$ is affine on each of the intervals $[a,b], [b,c],[c,d], [d,e]$,
\end{itemize}
where $[a,e]$ is the branch of $T$ with $a\in\IR$ and $a<b<c<d<e$. 
See Figure~\ref{fig:ex}
for a picture of the map $F$. This entirely determines $F$ because it is
of degree 1.

\begin{figure}[htb]
\centerline{\includegraphics{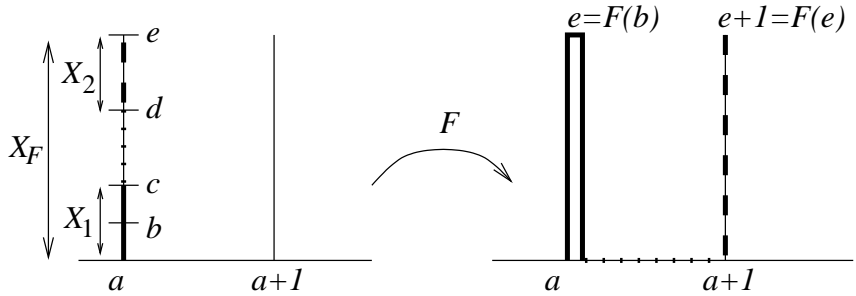}}
\caption{The action of $F$ on the branch $X_F$. \label{fig:ex}}
\end{figure}

$X_F$ is equal to $[a,e]$ and the intervals given by Lemma~\ref{lem:Xi}
are $X_1=[a,c]$ (with $p_1=0$) and $X_2=[d,e]$ (with $p_2=1$). $F$ is
an affine Markov map and the restriction of its Markov graph to $X_1,X_2$
is given in Figure~\ref{fig:Mgraph}. See \cite[Section 6.1]{AR} for general
results on Markov maps in this context, and in particular how it is possible
to deduce periodic {\modi} points and rotation numbers from the Markov graph.

\begin{figure}[htb]
\centerline{\includegraphics{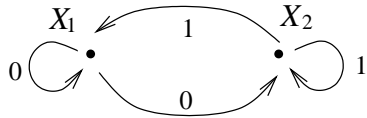}}
\caption{The Markov graph of $F$ restricted to the vertices $X_1$ and $X_2$
(actually, $X_1$ represents the two vertices $[a,b]$ and $[b,c]$). An
arrow $A\stackrel{i}{\longrightarrow}B$ means than $F(A)\supset B+i$.
\label{fig:Mgraph}}
\end{figure}

It can easily be deduced from the Markov graph of $F$ that $\Rot_{X_1\cap
X_{\infty}}(F)=[0,1]$ and the unique periodic {\modi} point $x\in X_F$
such that $\rho_F(x)=1$ is $x=e$, which does not belong to~$X_1$.

In addition, we note that $\Rot_{\IR}(F)=\{0\}$. Thus $\Rot(F)=[0,1]$ and $\Rot_{\IR}(F)$
is not a connected component of $\Rot(F)$.\hfill$\Box$
\end{ex}

\subsection{Rotation set of $F$}

Now, we are ready to prove that the set $\Rot(F)$ is closed and has finitely
many connected components, and that every rational number in $\Rot(F)$
is the rotation number of some periodic {\modi} point.
Note that in the following theorem, the intervals $I_0,\ldots, I_k$
may be not disjoint; in particular $I_0=\RotR(F)$ may not be a connected
component of $\Rot(F)$ (see Example~\ref{ex:counterex}).

\begin{theo}\label{theo:RotF}
Let $T\in\InfG$, $F\in\sigmap$ and 
$X_{\infty}=\{x\in X\mid \forall n\geq 1, F^n(x)\in X+\IZ\}$.
Then there exist an integer $k\geq 0$ and compact non empty intervals 
$I_0,\ldots, I_k$ in $\IR$ such that: 
\begin{itemize}
\item $I_0=\Rot_{\IR}(F)=\Rot_{\TR}(F)$,
\item
$\Rot(F)=\Rot^+(F)=\Rot^-(F)=I_0\cup\cdots\cup I_k$,
\item $\forall\; 1\leq i\leq k, \forall r\in I_i\cap\IQ$, there exists a periodic {\modi}
point $x\in  X_{\infty}$ with $\rhos(x)=r$,
\item $\forall\; 1\leq i\leq k,\ I_i\cap\IZ\not=\emptyset$.
\end{itemize}
Moreover, if $N$ is the integer given by Lemma~\ref{lem:Xi}, then $k\leq N$.
\end{theo}

\begin{demo}
Consider $X_1,\ldots, X_N$ and $p_1,\ldots,p_N$ given by Lemma~\ref{lem:Xi}. 
We define inductively $T_1,\ldots, T_N\subset T$ and $I_1,\ldots, I_N\subset\IR$ such that:
\begin{itemize}
\item[(a)]
$T_i$ is a connected subset of $T$ such that 
$T_i+1=T_i$, $F(T_i)\subset T_i$ and, for all $2\leq i\leq N$, $T_{i-1}\cup X_{i-1}
\subset T_i$.
\item[(b)] 
either $I_i$ is  empty, or $I_i$ is a compact interval containing
$p_i$ such that, for all $r\in I_i\cap\IQ$, there exists a 
periodic (mod 1) point $x\in  X_{\infty}$ with $\rhos(x)=r$.
\item[(c)] $\Rot_{T_i\cup X_i}(F)=\Rot_{T_i}(F)\cup I_i$,
and the same equality is valid with $\Rot^+$ and $\Rot^-$.
\item[(d)] If $i\geq 2$, 
$\Rot_{T_i}(F)=\Rot_{T_{i-1}\cup X_{i-1}}(F)$,
and the same equality is valid with $\Rot^+$ and $\Rot^-$.
\end{itemize}

Let $T_1=\TR$. It satisfies (a). If $X_1\cap X_{\infty}=
\emptyset$, we take $I_1=\emptyset$. Otherwise, we apply 
Proposition~\ref{prop:key} with $T'=T_1$, $Y=X_F$ 
and $X_1,\ldots, X_N$ is place
of $Y_1,\ldots,Y_M$. It provides a compact interval $I_1=I=
\Rot_{X_1\cap X_{\infty}}$ that satisfies (b). Moreover,
$\Rot_{T_1\cup X_1}(F)=\Rot_{T_1}(F)\cup \Rot_{X_1\cap X_{\infty}}(F)$,
and the same equality is valid with $\Rot^+$ and $\Rot^-$. Hence (c)
is satisfied for $i=1$.

\medskip
Let $i\geq 2$.
Suppose that $T_j$ and $I_j$ are already defined for all $1\leq j\leq i-1$, and
satisfy (a)-(d). Define
$$
A_i=[\min X_F,\min X_{i-1})\cup X_{i-1}
\cup\left(\left(\bigcup_{n\geq 1} F^n(X_{i-1})+\IZ\right)\cap X_F\right).
$$
For all $n\geq 1$, $F^n(\min X_{i-1})\in\TR$, and 
thus $(F^n(X_{i-1})+\IZ)\cap X_F$
is, either empty, or a compact subinterval of $X_F$ containing $\min X_F$.
Therefore, $A_i$ is a subinterval of $X$ containing $\min X_F$ and $X_{i-1}$.
Let $T_i=T_{i-1}\cup(A_i+\IZ)$. It is a connected subset of $T$, $T_i
+1=T_i$ and $T_{i-1}\cup X_{i-1}\subset T_i$.
Let us show that $F(T_i)\subset T_i$. Let $x\in T_{i-1}\cup A_i$. 
We distinguish 3  cases.
\begin{itemize}
\item If $x\in T_{i-1}$ then $F(x)\in T_{i-1}$ by invariance of $T_{i-1}$.
\item If $x\in [\min X_F, \min X_{i-1})$ then, either
$x\in X_1\cup\cdots\cup X_{i-2}\subset T_{i-1}$ and $F(x)\in T_{i-1}$,
or $x\in X\setminus(X_1\cup\cdots\cup X_N)$ and $F(x)\in \TR$.
\item If $x\in \left(\bigcup_{n\geq 0} F^n(X_{i-1})+\IZ\right)\cap X_F$
then, either $F(x)\in \left(\bigcup_{n\geq 0} F^n(X_{i-1})+\IZ\right)\cap (X_F+\IZ)\subset A_i+\IZ$, or $F(x)\in \TR$.
\end{itemize}
Consequently, $F(T_i)\subset T_i$, and (a) is satisfied. Moreover, 
$$
\Rot_{T_i}(F)=\Rot_{T_{i-1}}(F)\cup \Rot_{X_{i-1}}(F)=
\Rot_{T_{i-1}\cup X_{i-1}}(F),
$$
and the same equality is valid with $\Rot^+$ and $\Rot^-$, 
which is (d) for $i$.

If $F(X_i)\subset T_i$, we take $I_i=\emptyset$ and (b)-(c) are clearly
satisfied.
Otherwise, let $b\in X_i$ such that $F(b)\not\in T_i$. 
Let $Y=\overline{T\setminus T_i}\cap r^{-1}([0,1))$ and
$Y_{\infty}=\{x\in Y\mid \forall n\geq 1, F^n(x)\in Y+\IZ\}$.
The set $Y$ is
a compact subinterval of $X_F$ and $Y\cap\overline{T_i}=\{\min Y\}$. 
We have $b\in Y$ by invariance of $T_i$, and $F(b)\in Y+p_i$
because $F(X_i)\subset (X_F+p_i)\cup \TR$.

Let $a=\max (\min X_i,\min Y)$. 
We can define $c=\min\{x\in [a,b]\mid F(x)\in Y+\IZ\}$ because
$b\geq a$.
Moreover $F(a)\in \overline{T_i}$ because
$F(\min X_i)\in \TR$ and $F(\min Y)\in F(\overline{T_i})\subset 
\overline{T_i}$. 
Therefore, $F(c)=\min Y\ {\modi}$ by minimality.
Let $X_i'=[c,\max X_i]\subset X_i$. We can apply Proposition~\ref{prop:key}
with $T'=T_i$ and $X_i',X_{i+1},\ldots, X_N$ in place of
$Y_1,\ldots, Y_M$. We obtain a compact interval 
$I_i=I=\Rot_{X_i'\cap Y_{\infty}}(F)$ that satisfies (b) for $i$.

We have 
$$
X_i'\cap Y_{\infty}=\{x\in X_i\mid \forall n\geq 0, F^n(x)\not\in \Int{T_i}\}.
$$
Therefore, $\Rot_{T_i\cup X_i}(F)=\Rot_{T_i}(F)\cup 
\Rot_{X_i'\cap Y_{\infty}}(F)=\Rot_{T_i}(F)\cup I_i$,
and the same equality is valid with $\Rot^+$ and $\Rot^-$. Hence (c)
is satisfied for $i$. This concludes the construction of $T_i$ and $I_i$.

\medskip
Now, we end the proof of the theorem.
Since $F(T\setminus (T_N\cup X_N))\subset T_N$, it is clear that $\Rot(F)=
\Rot_{T_N\cup X_N}(F)$. Combining this with (c) and (d), we get that
\begin{eqnarray*}
\Rot(F)&=&\Rot_{T_N\cup X_N}(F)\\
&=&\Rot_{T_N}(F)\cup I_N\\
&=&\Rot_{T_{N-1}\cup X_{N-1}}(F)\cup I_N\\
&=&\Rot_{T_{N-1}}(F)\cup I_{N-1}\cup I_N\\
&=&\cdots\\
&=&\Rot_{T_1}(F)\cup I_1\cup\cdots\cup I_N
\end{eqnarray*}
and the same equalities are valid with $\Rot^+$ and $\Rot^-$.
Let $I_0=\Rot_{\IR}(F)$.
By Theorem~\ref{theo:rotR}, $I_0$ is a non empty compact interval
and $I_0=\Rot_{\TR}(F)=\Rot^+_{\TR}(F)=\Rot^-_{\IR}(F)$.
If $x\in T$, then either $x\in (X_1\cup\cdots\cup X_N)+\IZ$, or $F(x)\in T_1$.
Therefore, $\Rot_{T}(F)=\Rot_{T_1\cup X_1\cup\cdots\cup X_N}(F)$,
and the same equality is valid with $\Rot^+$ and $\Rot^-$.
To conclude, it remains to remove the empty intervals
among $I_1,\ldots, I_N$.
\end{demo}

If the empty rotation intervals are not removed in the proof of
Theorem~\ref{theo:RotF}, then the theorem can be stated as follows:

\medskip
\noindent\textbf{Theorem~\ref{theo:RotF}'}
Let $T\in\InfG$ and $F\in\sigmap$. Let
$X_1,\ldots,X_N$, $p_1,\ldots,p_N$ be given by Lemma~\ref{lem:Xi} and
$X_{\infty}=\{x\in X\mid \forall n\geq 1, F^n(x)\in X+\IZ\}$.
Then there exist compact intervals 
$I_0,\ldots, I_N$ in $\IR$ such that: 
\begin{itemize}
\item $I_0=\Rot_{\IR}(F)=\Rot_{\TR}(F)$,
\item
$\Rot(F)=\Rot^+(F)=\Rot^-(F)=I_0\cup\cdots\cup I_k$,
\item for all $1\leq i\leq N$, either $I_i=\emptyset$ or $p_i\in I_i$,
\item $\forall\; 1\leq i\leq N, \forall r\in I_i\cap\IQ$, there exists a 
periodic {\modi}
point $x\in \overline{\bigcup_{n\geq 0}(F^n(X_1)+\IZ)}\cap
X_{\infty}$ with $\rhos(x)=r$; if in addition
$r\in \Int{I_i}$ then $x$ can be chosen in $X_1\cap X_{\infty}$.
\end{itemize}

\bigskip
\noindent
Laboratoire de Math\'ematiques,\\
CNRS UMR 8628, B\^atiment 425,\\
Universit\'e Paris-Sud 11,\\
91405 Orsay cedex, France

\noindent
http://www.math.u-psud.fr/{\tiny$\sim$}ruette/\\
\texttt{Sylvie.Ruette@math.u-psud.fr}

\end{document}